\documentclass[12pt]{amsart}
\usepackage{latexsym, amsmath, amscd, amssymb, amsthm} 
\usepackage[english]{babel}
\usepackage{pb-diagram}
\newtheorem{te}{Theorem}[section]

 \newtheorem{lm}{Lema}[section]

\begin{document}

\noindent

 \title[ II. The symbolic method]{  The Weitzenb\"ock derivations and  classical invariant theory. II. The symbolic method }

\author{L. Bedratyuk}\address{Khmelnitskiy national university, Instituts'ka, 11,  Khmelnitskiy, 29016, Ukraine}

\begin{abstract} A  method based on the  symbolic methods of the classical invariant theory  is developed for a  representation of elements of kernel of Weitzenb\"ok derivations.
\end{abstract}
\maketitle

\section{Introduction}
 Let  $\mathbb{K}$ be a field of characteristic 0 and let   $\mathbb{K}[X]$ be  a polynomial algebra in a set of variables  $X.$  A linear locally
nilpotent derivation $\mathcal{D}$
of the polynomial algebra $\mathbb{K}[X]$ 
is called a Weitzenb\"ock derivation.   
 Denote by  $\mathcal{D}_{\mathit{\mathbf{d}}},$  $\mathit{\mathbf{d}}:=(d_1,d_2,\ldots, d_s)$   the   Weitzenb\"ok derivation of the algebra $\mathbb{K}[X]$ with   the Jordan normal form  consisting  of  $s$  Jordan blocks of size  $d_1+1,$ $d_2+1,$ $\ldots, d_s+1,$ respectively.  The only   derivation which corresponds to  a single  Jordan block of size $d+1$ is called the basic Weitzenb\"ock derivation and denoted by $\mathcal{D}_d.$  
The algebra  
$$
{\rm{ker}}\, \mathcal{D}_{\mathit{\mathbf{d}}}=\left\{ f \in \mathbb{K}[X] \mid  \mathcal{D}_{\mathit{\mathbf{d}}}(f)=0 \right\},
$$
is called {\it the kernel}  of the derivation  $\mathcal{D}_{\mathit{\mathbf{d}}}.$ 
It is well known that the kernel   $\ker \mathcal{D}_{\mathit{\mathbf{d}}}$  is a finitely generated algebra, see \cite{Wei}--\cite{Tyc}.  
  However, it remained an open problem to find a minimal system of homogeneous
generators (or even the cardinality of such a system) of the algebra $\rm{ker}\,D_{\mathit{\mathbf{d}}}$ even  for small tuples  ${\mathit{\mathbf{d}}}.$ 

The aim of this paper is to develop an effective method for the  representation and manipulation of the kernel elements of Weitzenb\"ok derivations.

 In a previous work   \cite{CIT-LND_I} we   showed that the kernel of the derivations $\mathcal{D}_{\mathit{\mathbf{d}}},$  $\mathit{\mathbf{d}}:=(d_1,d_2,\ldots, d_s)$ is isomorphic to the algebra joint covariants of   $s$ binary forms of orders  $d_1,d_2,\ldots, d_s.$  Algebras of joint covariants of  binary forms   were an  object of research in the invariant theory
in  the 19th century. To describe the kernel of linear locally  nilpotent derivations we  should  involve a computational tools of classical invariant theory, including the famous symbolic method. The symbolic method  was developed by Aronhold, Clebsch, and Gordan.  It is    is a most powerful tool of the classical invariant theory. A classic presentation of the symbolic method can be found in
\cite{Cleb}--\cite{Gle}.
Recently, a rigorous foundation of the symbolic method   has been given by Kung and Rota \cite{KR} and by Kraft and Weyman \cite{KW}.

In this paper we develop an analogue of the classical symbolic method for the kernel of Weitzenb\"ok derivations.
In order to explain the essence of the method, we give some examples.
 Let   $\mathcal{D}_{4}$ be  the basic Weitzenb\"ock derivation of the polynomial algebra   ${\mathbb{K}[X_4]:=\mathbb{K}[x_0,x_1,x_2,x_3,x_4]}$ i.e., $\mathcal{D}_{4}(x_i)=x_{i-1},$ $\mathcal{D}_{4}(x_0)=0,$ $i=0\ldots 4.$  Consider the Weitzenb\"ock derivation  $\mathcal{D}_{(4,4)}$ of  the algebra $\mathbb{K}[X_4,Y_4]:$ $\mathcal{D}_{(4,4)}(x_i)= \mathcal{D}_{4}(x_i),$  $\mathcal{D}_{(4,4)}(y_i)=y_{i-1},$ $\mathcal{D}_{(4,4)}(y_0)=0.$ 
The following differential operator  $\mathcal{P}_{y,x}^{(n)}:\mathbb{K}[X_n]  \to \mathbb{K}[X_n,Y_n]$ defined by
$$
\mathcal{P}_{y,x}^{(n)}=y_0 \frac{\partial}{\partial  x_0}+y_1 \frac{\partial}{\partial  x_1}+\cdots +y_n \frac{\partial}{\partial  x_n},
$$
 is called the {\it polarization operator.}  The operator  $\mathcal{R}_{y,x}^{(n)}:\mathbb{K}[X_n,Y_n]  \to \mathbb{K}[X_n],$   defined by $${\mathcal{R}_{y,x}^{(n)}(F)=F\Bigl |_{y_i=x_i}}$$ is called the  {\it restitution operator.} If  $F$ is a homogeneous polynomial then Euler's homoge\-neous function theorem implies that  $\mathcal{R}_{y,x}^{(n)}\left(\mathcal{P}_x^y(F)\right) =\deg(F) \,F.$ The polarization operator commutes with the Weitzenb\"ock derivations:
$$
\mathcal{P}_{y,x}^{(4)}\left( \mathcal{D}_4(F)\right)=\mathcal{D}_{(4,4)}\left(\mathcal{P}_{y,x}^{(4)}(F) \right).
$$

It is easy to verify that the polynomial  $F={x_{{2}}}^{2}+2\,x_{{0}}x_{{4}}-2\,x_{{1}}x_{{3}}$ belongs to $\ker \mathcal{D}_{4}$ and its polarization  $$\mathcal{P}_{y,x}^{(4)}(F)=2\,y_{{0}}x_{{4}}-2\,y_{{1}}x_{{3}}+2\,y_{{2}}x_{{2}}-2\,y_{{3}}x_{{1}
}+2\,y_{{4}}x_{{0}},$$
belongs to  $\ker \mathcal{D}_{(4,4)}.$ Let us  change the   variables  by
\begin{equation}\label{ch_v}
x_i=\frac{1}{i!} \alpha_0^{4-i} \alpha_1^i,y_i=\frac{1}{i!} \beta_0^{4-i} \beta_1^i, i=0,\ldots,4.
\end{equation}
Then we get  that 
$$
\mathcal{P}_{y,x}^{(4)}(F)=\frac{1}{12} \,  \left( \alpha_{{0}}\beta_{{1}}-\beta_{{0}}\alpha_{{1}} \right) := \frac{1}{12}[\alpha, \beta]^4,
$$
where $[\alpha, \beta]:= \alpha_{{0}}\beta_{{1}}-\beta_{{0}}\alpha_{{1}}.$
The  polynomial  $\Psi=\frac{1}{12}[\alpha, \beta]^4$ is called  the  {\it  symbolic representation} of the polynomial  $F.$ The letters  $\alpha, \beta$ are called the {\it  symbol letters.} Observe, that $\Psi$  belongs to kernel of the  derivation $\mathcal{D}_{(1,1)}$ which acts on  $\mathbb{K}[\alpha_0,\alpha_1,\beta_0,\beta_1].$ Moreover, the  polynomial $\Psi$ has much simpler form than the polynomial  $F.$

On the other hand, let us consider the polynomial   $\Phi=\alpha_0^2 \beta_0^2 [\alpha, \beta]^2 \in \ker \mathcal{D}_{(1,1)}.$ Then by (\ref{ch_v})  we get 
$$
\Phi={\alpha_{{0}}}^{4}{\beta_{{0}}}^{2}{\beta_{{1}}}^{2}-2\,{\alpha_{{0}}}
^{3}{\beta_{{0}}}^{3}\beta_{{1}}\alpha_{{1}}+{\alpha_{{0}}}^{2}{\beta_
{{0}}}^{4}{\alpha_{{1}}}^{2}=2(x_0 y_2-x_1 y_1+x_2 y_0),
$$
and  $\mathcal{R}_{y,x}^{(4)}(\Phi)=2\,(2 x_0 x_2-x_1^2)\in \ker \mathcal{D}_{4}.$  Thus  $ \frac{1}{2} \alpha_0^2 \beta_0^2 [\alpha, \beta]^2$ is a symbolic representation   for  $2 x_0 x_2-x_1^2.$  To get elements of degree 3 we  should involve  one  more symbolic letter  $\gamma.$ Similarly, one may show that   $${\frac{1}{4!}[\alpha, \beta]^2 [\alpha,\gamma]^2 [\beta,\gamma]^2 \in \ker \mathcal{D}_{(1,1,1)}},$$  is a symbolic representation for the following element of the kernel of derivation $\mathcal{D}_{4}:$
$$
12\,x_{{0}}x_{{4}}x_{
{2}}+6\,x_{{1}}x_{{3}}x_{{2}}-2\,{x_{{2}}}^{3}-9\,x_{{0}}{x_{{3}}}^{2}-6\,{x_{{1}}}^{2}x_{{4}},
$$
and for all its polarizations.

For the general case consider the  polynomial algebra $\mathbb{K}[\alpha_i \mid \alpha \in \mathcal{J}, i=0,1]$  where $\alpha$ runs a set $\mathcal{J}$  of symbol letters.  Elements of kernel of the Weitzenb\"ock derivation  $$\mathcal{D}_{\mathcal{J}}:=\mathcal{D}_{(\underbrace{1,1,\ldots,1}_{ |\mathcal{J}|\, {\rm   times} })}$$ defined by $\mathcal{D}_{\mathcal{J}}(\alpha_0)=0,$ $\mathcal{D}_{\mathcal{J}}(\alpha_1)=\alpha_0,$ for all $\alpha \in \mathcal{J}$ are  called  the symbolic expressions.

The following statement is a main point of the symbolic method:

{{\bf The Symbolic Method.} {\it  
\begin{itemize}
	\item Any element of ${\rm{ker}}\, \mathcal{D}_{\mathit{\mathbf{d}}}$  allows a symbolic representation;

 \item Any symbolic  expression is a symbolic representation for an  element of ${\rm{ker}}\, \mathcal{D}_{\mathit{\mathbf{d}}}$  for some $\mathit{\mathbf{d}}.$
\end{itemize}
}}
Thus  we get  a remarkable fact --  the   kernel of an arbitrary Weitzenb\"ock derivation $\mathcal{D}_{\mathit{\mathbf{d}}}$ is completely defined by the kernel of the special Weitzenb\"ock derivation $\mathcal{D}_{\mathcal{J}}.$ 
The kernel  $ \ker \mathcal{D}_{\mathcal{J}}$ is well-known and  generated  by  $\alpha_0$ and the brackets $[\alpha,\beta]$  where $\alpha, \beta $  run over $\mathcal{J}.$

 The paper organized as follows. In section 2 we review some of the standard 
facts on the representation theory of the Lie algebra  $\mathfrak{sl}_2$ and its maximal nilpotent subalgebra $\mathfrak{u}_2.$ 
  
  In section 3 we develop an analogue of the classical symbolic method for Weitzenb\"ock derivations.   In section 4  we  introduce  the notions of the convolution and the semi-transvectant  which used for calculation a generating set of kernel of derivations.

%%%%%%%%%%%%%%%%%%%%%%%%%%%%%%%%%%%%%%%%%%%%5
\section{Basic facts}

A  representation  of the Lie algebra  $\mathfrak{g}$ on  a  finite-dimensional complex vector space  $V$  is a homomorphism  $\rho : \mathfrak{g} \to \mathfrak{gl}(V),$  where  $\mathfrak{gl}(V)$ is Lie algebra of endomorphisms of $V.$  We say that such a map gives $V$ the structure of  $\mathfrak{g}$-module.The algebra   $\mathfrak{g}$  acts on $V$  by linear operators  $\rho(g), g \in \mathfrak{g}.$  When 
there is little ambiguity about the map $\rho$  we sometimes call $V$  itself a representation of $\mathfrak{g}$; in this vein we will  suppress the symbol $\rho$ and write  $g\,v$ for $\rho(g)v.$ 

If  $U,V$ are representations,  the tensor product  $U \otimes V$  is also representation, the latter via   
 $$g( u \otimes v)=gu \otimes v+u \otimes gv.$$
 For a representation $V,$ the  tensor algebra ${\rm T}(V)$ is again a representation of 
$\mathfrak{g}$ by this rule,  and symmetric algebra  ${\rm Sym}(V)$ 
are subrepresentations of it. 
Thus, the algebra   $\mathfrak{g}$ acts on  ${\rm Sym}(V)$  by derivations. An element  $v \in V$  is called an  {\it invariant }  of   $\mathfrak{g}$-module  $V$ if  $gv=0.$ Denote by  $V^{\mathfrak{g}}$  the set of all invariants of the $\mathfrak{g}$-module $V.$

Let  $\mathfrak{sl}_2$ is the Lie algebra of  $2\times 2$ traceless matrices and  $\mathfrak{u}_2$ is its maximal nilpotent subalgebra.  The canonical basis of  $\mathfrak{sl}_2$  is the basis $(e,f,h),$  where
$$
e=\begin{pmatrix} 0 & 1\\ 0 & 0 \end{pmatrix}, f=\begin{pmatrix} 0 & 0\\ 1 & 0 \end{pmatrix}, h=\begin{pmatrix} 1 & 0\\ 0 & -1 \end{pmatrix}.
$$
We  have   
\begin{equation}\label{sl2}
[h,e]=2\,e, [h,f]=-2\,f,[e,f]=h.
\end{equation}
For each nonnegative integer $n,$ the algebra  $\mathfrak{sl}_2$ has an irreducible representation $V_n$ of dimension $n+1,$ which is unique up to an isomorphism.
 The endomorphisms  $ \rho(e),$ $\rho(f),$ $\rho(h)$ act on  $V_n:=\left\langle v_0,v_1,\ldots,v_n \right\rangle$ by the formulas
 $$
 \rho(e)v_i= v_{i-1},  \rho(e)v_i=(n-i)(i+1) v_{i+1}, \rho(h)v_i=(n-2i) v_i.
 $$
The action of $\mathfrak{sl}_2$ extends by derivations to  the symmetric algebra  ${\rm Sym}(V_n)$ and to  the algebra  $${\rm Sym}(V_{\bf d}):={\rm Sym}(V_{d_1} \oplus V_{d_2}\oplus \ldots \oplus V_{d_s}),\mathit{\mathbf{d}}:=(d_1,d_2,\ldots, d_s).$$
For convenience,  the derivations of ${\rm Sym}(V_{\bf d})$ which  correspond to the operators 
 $ \rho(e),$ $\rho(f),$ $\rho(h)$  denote by  ${\bf D},$ ${\bf D}_*$ and ${\bf E}$  respectively.
 Let us identify the algebra ${\rm Sym}(V_n)$  with the polynomial algebra   $\mathbb{K}[V_n]:=\mathbb{K}[v_0,v_1,\ldots,v_n]$  and  identify the algebra  ${\rm Sym}(V_{\bf d})$ with the polynomial algebra $\mathbb{K}[V_{d_1},V_{d_2},\ldots V_{d_s}]$  Under this identification, the kernel of the derivation  ${\bf D}$  coincides with the algebra of invariants ${\rm Sym}(V_{\bf d})^{\mathfrak{u}_2},$ and the algebra  $\ker {\bf D} \cap \ker {\bf D}_*$  coincides with the algebra of invariants  ${\rm Sym}(V_{\bf d})^{\mathfrak{sl}_2}.$  Any   $\mathfrak{u}_2$-invariant  is called the   {\it semi-invariant}. It is clear that the derivation   ${\bf D}$  is exactly  the Weitzenb\"ock derivation $\mathcal{D}_{\mathit{\mathbf{d}}}.$
 
To begin, let us describe the algebras of invariants and  semi-invariants of the symmetric algebra  ${\rm Sym}(V_{1} \oplus V_{1}\oplus \ldots \oplus V_{1}).$
 
Let  $\mathfrak{G} = \{\alpha, \beta, \ldots  \}$  be  an alphabet consisting of an infinite supply of Greek letters. The letter in  $\mathfrak{G}$ are called the {\it symbol letter}. To each symbol letter 
 $\alpha$ we associate two variables   $\alpha_0,$  $\alpha_1$ and the two-dimensional vector space  $V_{\alpha}:=\mathbb{K} \alpha_0 \oplus \mathbb{K} \alpha_1.$  For a finite subset  $\mathcal{J} \subset \mathfrak{G}$ put  $V_{\mathcal{J}}:=\oplus_{\alpha \in J} V_{\alpha}.$ The algebra  ${\rm Sym}\left(V_{\mathcal{J}} \right)$ turns into  $\mathfrak{sl}_2$-module by  the actions:
\begin{equation}\label{rul_1}
{\bf D}(\alpha_0)=0, {\bf D}(\alpha_1)=\alpha_0, {\bf D}_*(\alpha_0)= \alpha_{1},{\bf D}_*(\alpha_2)= 0, {\bf E}(\alpha_i)=(1-2i) \alpha_i, i=0,1, \alpha \in \mathcal{J}.
\end{equation}
 A direct check shows  that the conditions  (\ref{sl2}) hold.
We  note that the symmetric group $S_{ |\mathcal{J}|}$ naturally  acts  on ${\rm Sym}\left(V_{\mathcal{J}} \right),$ here $|\mathcal{J}|$ is the cardinality of the set $\mathcal{J}.$

{\bf The First Fundamental Theorem for  $\mathfrak{sl}_2.$} {\it The algebra of semi-invariants  ${\rm Sym}\left(V_{\mathcal{J}} \right)^{\mathfrak{u}_2}$  is generated by  $\alpha_0$ and by the brackets  $$[\alpha, \beta]:=\begin{vmatrix}
\alpha_0 & \alpha_1 \\
\beta_0 & \beta_1
\end{vmatrix}, \alpha, \beta \in \mathcal{J}. $$  The algebra  of invariants ${\rm Sym}\left(V_{\mathcal{J}} \right)^{\mathfrak{sl}_2}$ is generated by the brackets $[\alpha, \beta].$}

It is  well-known results of classical invariant theory, see  \cite{GrY}. In the theory of locally nilpotent derivations the  first part of this statement is known as the Nowicki conjecture,  see,  for instance,  \cite{BNS}.

{\bf Corollary.} {\it As a vector space the algebra  ${\rm Sym}\left(V_{\mathcal{J}} \right)^{\mathfrak{u}_2}$ is generated by the polynomials
$$
P_{l,m}:=\prod_{\alpha \neq \beta}[\alpha, \beta]^{l_{\alpha, \beta}} \prod_{\gamma} \gamma_0^{m_{\gamma}}, \alpha, \beta, \gamma \in \mathcal{J}.
$$
}
 The {\it order}  ${\rm ord}\, P$  and the {\it weight} ${\rm  wt}\, P  := ({\rm wt}_{\alpha} P)_{\alpha \in \mathcal{J}}$  of the
symbolic expression $P$  are defined by 
$$
{\rm ord}\, P:=\sum_{\gamma} m_{\gamma}, {\rm wt}_{\alpha} \, P:=\sum_{\beta}\left(l_{\alpha,\beta}+l_{\beta,\alpha}\right)+m_{\alpha}.
$$
In particular, ${\rm wt}_{\alpha} P$ is equal to the number of times the symbol $\alpha$  occurs in the
symbolic expression $P.$  Observe, that 
$$
{\rm ord}\, P=\min_k\{k \in \mathbb{N} \mid {\bf D}_*^{k+1}(P)=0\},\, {\bf E}\left(P\right)={\rm  wt}\, P\,\cdot  P.
$$

The symbolic expression $P$ is called {\it decomposable} if it can be written as a product
$P = P_1\,P_2$ in a non-trivial way where $P_1$ and $P_2$ are {\it disjoint,} i.e., no symbol occurs
in both. We denote by ${\rm supp}\, P$ the support of $P,$ i.e., the set of symbols
$\alpha \in \mathcal{J}$  occurring in $P.$  Put  $\alpha \sim \beta$ if  ${\rm wt}_{\alpha}={\rm wt}_{\beta}.$  Then the relation $\sim$ is an equivalence relation  defined on the set ${\rm supp}\, P \subseteq \mathcal{J}.$ Denote by  $\mathcal{J}_1,$ $\mathcal{J}_2,\ldots,$ $\mathcal{J}_t$ the equivalence classes and by   $m_1,$ $m_2,$  $\ldots m_t$ denote their cardinality, $t=|\,{\rm supp}\, P/\sim|.$ Denote by  $n_1,$ $n_2,$  $\ldots n_t$  the corresponding weights of elements of the classes $J_i.$

Let $x_{\alpha, \beta}, x_{\gamma}$ for $\alpha, \beta, \gamma \in \mathcal{J},$ $\alpha \neq \beta$ denote independent variables and define the free polynomial algebra 
$$
{\rm Sym}_{\mathcal{J}}:=\mathbb{K}[x_{\alpha, \beta}, x_{\gamma} \mid \alpha, \beta, \gamma \in \mathcal{J},\alpha \neq \beta].
$$ 
Define  the map $\chi:{\rm Sym}_{\mathcal{J}} \to {\rm Sym}\left(V_{\mathcal{J}} \right)^{\mathfrak{u}_2}$ by 
$$
\chi(x_{\alpha, \beta})=[\alpha, \beta], \chi(x_{\gamma})=\gamma_0,
$$
and extend it in the natural way to monomials and all of ${\rm Sym}_{\mathcal{J}}.$

{\bf The Second Fundamental Theorem $\mathfrak{sl}_2.$} {\it There is a canonical isomorphism
$${\rm Sym}_{\mathcal{J}}/\ker \chi \cong {\rm Sym}\left(V_{\mathcal{J}} \right)^{\mathfrak{u}_2}$$

where the ideal $\ker \chi$  is generated by the elements
\begin{align*}
&x_{\alpha, \beta}+x_{\beta, \alpha}=0,x_{\gamma} x_{\alpha,\beta}+x_{\beta} x_{\gamma, \alpha} + x_{\alpha} x_{\beta, \gamma}=0,\\
&x_{\alpha, \beta} x_{\gamma, \delta}+x_{\gamma, \alpha} x_{\beta, \delta}+x_{\beta, \gamma} x_{\alpha, \delta}=0.
\end{align*}
}
The theorem implies  that the following three relations(syzygies)  
\begin{gather*}\label{syg}
[\alpha, \beta]+[\beta, \alpha]=0,  \\
\gamma_0 [\alpha,\beta]+\beta_0 [\gamma, \alpha]+\alpha_0 [\beta, \gamma]=0,\\
[\alpha, \beta][\gamma, \delta]+[\gamma, \alpha] [\beta, \delta]+[\beta, \gamma] [\alpha, \delta]=0,
\end{gather*}
for distinct $\alpha, \beta, \gamma, \delta \in  \mathcal{J} $ generate all the relationship among the semi-invariants ${\rm Sym}\left(V_{\mathcal{J}} \right)^{\mathfrak{u}_2}.$

The algebra ${\rm Sym}\left(V_{\mathcal{J}} \right)^{\mathfrak{u}_2}$ is graded  by weight. Let ${\rm Sym}\left(V_{\mathcal{J}} \right)^{\mathfrak{u}_2}_{{\bf w}}$  be   those elements with weight ${\bf w}.$   Then ${\rm Sym}\left(V_{\mathcal{J}} \right)^{\mathfrak{u}_2}=\oplus_{{\bf w}} {\rm Sym}\left(V_{\mathcal{J}} \right)^{\mathfrak{u}_2}_{{\bf w}}.$  If ${\bf w}=(n,n,\ldots, n)$  we will write ${\bf w}=(n)^{|\mathcal{J}|}.$   If 
$$
{\bf w}=(\underbrace{n_1,n_1,\ldots, n_1,}_{m_1 \, \text{times}} \underbrace{n_2,n_2,\ldots, n_2,}_{m_2 \, \text{times}} \ldots  \underbrace{n_t,n_t,\ldots, n_t}_{m_t \, \text{times}}),
$$  
we write ${\bf w}=(n_1)^{m_1} (n_2)^{m_2} \cdots (n_t)^{m_t}, $  or more compact ${\bf w}={\bf n}^{{\bf m}}.$ Here ${\bf n}:=(n_1, n_2, \ldots n_t),$ ${\bf m}:=(m_1, m_2, \ldots m_t)$ and $m_1+m_2+\cdots +m_t=|\mathcal{J}|.$

\section{Symbolic method}

%%%%%%%%%%%%%%%%%%%%%%%%%%%%%%%%%%%%%%%%%%%%%%%%
We will show that symbolic expressions can be used in a very efficient way to
describe and manipulate semi-invariants of the ${\mathfrak{u}_2}$-module ${\rm Sym}(V_{\bf d}),$ $\mathit{\mathbf{d}}:=(d_1,d_2,\ldots, d_s).$   Recall that  ${\rm Sym}(V_{\bf d})^{{\mathfrak{u}_2}}=\ker \mathcal{D}_{\mathit{\mathbf{d}}}.$

Let  $\mathfrak{R} = \{x, y,z,  \ldots \}$ be  an alphabet consisting of an infinite supply of ordered roman  letters.  Denote by  $n_x$ the ordinal number of the letter $x$  in the set  $\mathfrak{R}.$ To each letter $x$ and to each  integer number $n$ associate  
the  $n+1$-dimension vector space  $$V_{x,n}:=\mathbb{K} x_0 \oplus \mathbb{K} x_1 \oplus \cdots \oplus \mathbb{K}x_n\cong V_n.$$  For a finite subset  $\mathcal{I} \subset \mathfrak{R}$   and for   $\mathit{\mathbf{d}}:=(d_1,d_2,\ldots, d_{|\mathcal{I}|})$ put $V_{\mathcal{I},\mathit{\mathbf{d}}}:=\oplus_{x \in \mathcal{I}} V_{x,d_i}.$ The action  
\begin{equation}\label{rul_2}
{\bf D}(x_i)= x_{i-1},  {\bf D}_*(x_i)=(i+1)(d_{n_x}-i) x_{i+1}, {\bf E}(x_i)=(d_{n_x}-2i) x_i, i=0,1,\ldots, d_{n_x},  x \in \mathcal{I}.
\end{equation}
gives  ${\rm Sym}\left(V_{\mathcal{I,\mathit{\mathbf{d}}}} \right)$ the structure of $\mathfrak{sl}_2$-module.
The {\it order} ${\rm ord}\, S$   of a homogeneous semi-invariant $S \in {\rm Sym}\left(V_{\mathcal{I,\mathit{\mathbf{d}}}} \right) ^{{\mathfrak{u}_2}}$  is defined  by 
$$
{\rm ord}\,S:=\min_k \{ k \mid {\bf D}_*^{k+1}(S)=0 \}.
$$

The algebra ${\rm Sym}\left(V_{\mathcal{I,\mathit{\mathbf{d}}}} \right) ^{{\mathfrak{u}_2}}$ is graded  by multidegree. Let $\left({\rm Sym}\left(V_{\mathcal{I,\mathit{\mathbf{d}}}} \right) ^{{\mathfrak{u}_2}}\right)_{\mathit{\mathbf{m}}}$  be   those elements with multidegree $\mathit{\mathbf{m}}.$  Then 
$${\rm Sym}\left(V_{\mathcal{I,\mathit{\mathbf{d}}}} \right) ^{{\mathfrak{u}_2}}=\oplus_{\mathit{\mathbf{m}}} \left({\rm Sym}\left(V_{\mathcal{I,\mathit{\mathbf{d}}}} \right) ^{{\mathfrak{u}_2}}\right)_{\mathit{\mathbf{m}}}.
$$  

 The following result  summarizes what is classically
called "symbolic method":

%+++++++++++++++++++++++++++++++++++++++++++++++
\begin{te}\label{main} There is a surjective ${\mathfrak{u}_2}$-homomorphism of vector spases
$$
\Lambda: \left({\rm Sym}\left(V_{\mathcal{J}} \right)^{{\mathfrak{u}_2}} \right)_{{\bf d}^{{\bf m}}} \to \left({\rm Sym}\left(V_{\mathcal{I,\mathit{\mathbf{d}}}} \right) ^{{\mathfrak{u}_2}}\right)_{\mathit{\mathbf{m}}},{\bf m}:=(m_1, m_2, \ldots m_t), |\mathcal{I}|=t,
$$
such that the composition $\Lambda \circ \chi:\left({\rm  Sym}_{\mathcal{J}} \right)_{{\bf d}^{{\bf m}}} \to \left({\rm Sym}\left(V_{\mathcal{I,\mathit{\mathbf{d}}}} \right) ^{{\mathfrak{u}_2}}\right)_{\mathit{\mathbf{m}}}$  is surjective with kernel 
$$
\ker \Lambda \circ \chi= \left(\ker \chi \right)_{{\bf d}^{{\bf m}}}+\left\langle P-\sigma P \mid P \in \left({\rm  Sym}_{\mathcal{J}} \right)_{{\bf d}^{{\bf m}}} , \sigma \in S_{|\mathcal{J}_i|}, \mathcal{J}_i \in {\rm Supp}\,P/\sim, i=1,\ldots t \right\rangle.
$$
\end{te}
%+++++++++++++++++++++++++++++++++++++++++++++++++++

 \begin{proof}
 Let $V_{\alpha}=\left\langle \alpha_0, \alpha_1 \right\rangle $ and $V_{x,n}=\left\langle x_0, x_1, \ldots x_n \right \rangle$ be  two ${\mathfrak{sl}_2}$-modules  as above.  It is well-known that the linear map $\alpha_0^{n-i} \alpha_1^i  \longmapsto i! x_{i}$
 is  ${\mathfrak{sl}_2}$-isomorphism of the vector spaces   ${\rm Sym}^n(V_{\alpha})$  and   $V_{x,n}.$
 In fact, we have
\begin{align*}
 &{\bf D}(\alpha_0^{n-i} \alpha_1^i)=i \alpha_0^{n-i+1} \alpha_1^{i-1} \longmapsto i (i-1)!\, x_{i-1}=i!\, {\bf D}(x_i),\\
 &{\bf D}_*(\alpha_0^{n-i} \alpha_1^i)=(n-i) \alpha_0^{n-i-1} \alpha_1^{i+1} \longmapsto (n-i) (i+1)!\, x_{i+1}=i! \,{\bf D}_*(x_i).
\end{align*}
 
 Let us consider a set  $\mathcal{J} \subset \mathfrak{G},$ $|\mathcal{J}|=n.$  The linear multiplicative map $\alpha_0^{d-i} \alpha_1^i  \longmapsto i!\,x_{i}$  for all $\alpha \in \mathcal{J} $  determines the ${\mathfrak{sl}_2}$-homomorphism  of the component  $\left({\rm Sym}\left(V_{\mathcal{J}} \right)\right)_{{\bf w}=(d)^n}$  into ${\rm Sym}^n(V_d).$
 
Let us now consider the component $\left({\rm Sym}\left(V_{\mathcal{J}} \right)^{{\mathfrak{u}_2}}\right)_{{\bf w}},$  where ${\bf w}=(n_1)^{m_1} (n_2)^{m_2} \cdots (n_t)^{m_t}$  and $m_1+\cdots+m_t=|\mathcal{J}|.$  Let $P$ be  a symbol expression of 
 $\left({\rm Sym}\left(V_{\mathcal{J}} \right)^{{\mathfrak{u}_2}}\right)_{{\bf w}}.$   Let $\varphi$ be a surjective map of  the coset ${\rm supp} P/\sim$ into a finite set $\mathcal{I} \subset \mathfrak{R}.$  Define the  map $\Lambda$ by  $$\alpha_0^{{\rm wt}_{\alpha}-i} \alpha_1^i  \mapsto i!\,\varphi(\alpha)_i, $$ $i=0,1,\ldots {\rm wt}_{\alpha}$ for all  $\alpha \in {\rm supp}\, P$ and extend it in the natural way to monomials.  Then  $\Lambda(P)$  is a multihomogeneous polynomial in the set of   $t$ letters   $\mathcal{I}$ of  the multidegree  ${\bf m}:=(m_1, m_2, \ldots m_t).$  To each  roman letter  $x \in \mathcal{I}$  we associate the variable set  $x_0, x_1, \ldots, x_{{\rm wt}_{\alpha}},$ where $\varphi(\alpha)=x,$ and ${\rm wt}_{\alpha}=n_i$  for  some $i.$  Therefore we can conclude that  $\Lambda(P)  \in \left({\rm Sym}\left(V_{\mathcal{I,\mathit{\mathbf{d}}}} \right) \right)_{\mathit{\mathbf{m}}}.$  Since  $\Lambda$ is  ${\mathfrak{u}_2}$-homomorphism and  ${\bf D} \left(P \right)=0,$  we see that the following inclusion holds: 
 $$
 \Lambda \left( \left({\rm Sym}\left(V_{\mathcal{J}} \right)^{{\mathfrak{u}_2}}\right)_{{\bf w}} \right) \subseteq \left({\rm Sym}\left(V_{\mathcal{I,\mathit{\mathbf{d}}}} \right) ^{{\mathfrak{u}_2}}\right)_{\mathit{\mathbf{m}}}.
 $$
Let us prove the surjectivity of the map $\Lambda.$   Rewrite the $\mathfrak{u}_2$-module $V_{\mathcal{I,\mathit{\mathbf{d}}}}$  as  $V_{\mathcal{I,\mathit{\mathbf{d}}}}=\oplus_{\alpha} m_i V_{{\varphi(\alpha)}, n_i},$  where $\alpha$  runs over  all members of cosets.    Every semi-invariant
 is a sum of multihomogeneous semi-invariants. Moreover, a multihomogeneous
semi-invariant $S \in \left({\rm Sym}\left(V_{\mathcal{I,\mathit{\mathbf{d}}}} \right) ^{{\mathfrak{u}_2}}\right)_{\mathit{\mathbf{m}}}$ of multidegree $(m_1, m_2, \ldots m_t)$ can be polarized   to produce a {\it multilinear} semi-invariant $\mathcal{P}(S)$ of multidegree $(\underbrace{1,1, \ldots, 1}_{|\mathcal{J}| \text{ times }}).$  Note that 
$$
\mathcal{P}(S) \in \left({\rm Sym}\left( \widetilde V_{\mathcal{I,\mathit{\mathbf{d}}}} \right) ^{{\mathfrak{u}_2}}\right)_{(1,1, \ldots, 1)},
$$
where  $\widetilde V_{\mathcal{I,\mathit{\mathbf{d}}}}=\oplus_{x \in \widetilde{\mathcal{I}}} V_{x,n_i}, $ $\widetilde{\mathcal{I}}$ is union $\mathcal{I}$  with the set of  new polarizing variables, $|\widetilde{\mathcal{I}}|=|\mathcal{J}|.$ Clearly, $\mathcal{P}(S)$ can be reconstructed from $S$  by  restitution.
Since $|\widetilde{\mathcal{I}}|=|\mathcal{J}|$  we can associate to each  symbol letter $\alpha \in \mathcal{J}$ a vector space $V_{x,n^*},$ for some natural $n^*.$ Extend the map  $x_k \mapsto 1/k!\,\alpha_0^{n^*-k} \alpha_1^k$ multiplicatively to all monomial and denote it by  $\widetilde{\Lambda}.$  We have the commutative diagram
\[ \begin{diagram}
\node{\left({\rm Sym}\left(V_{\mathcal{J}} \right)^{{\mathfrak{u}_2}} \right)_{{\bf d}^{{\bf m}}}} \arrow[2]{e,t}{\Lambda}
\node[2]{\left({\rm Sym}\left(V_{\mathcal{I,\mathit{\mathbf{d}}}} \right) ^{{\mathfrak{u}_2}}\right)_{\mathit{\mathbf{m}}}} 
\arrow{sw,b}{\mathcal{P}}
\\
\node[2]{\left({\rm Sym}\left( \widetilde V_{\mathcal{I,\mathit{\mathbf{d}}}} \right) ^{{\mathfrak{u}_2}}\right)_{(\underbrace{1,1, \ldots, 1}_{|\mathcal{I}|\text{ times} })}} \arrow{nw,r}{\widetilde{\Lambda}}
\end{diagram}\]
where all arrows are ${\mathfrak{sl}_2}$-homomorphisms, and  $\Lambda\left( \widetilde{\Lambda}(\mathcal{P}(S))\right)=S.$  Thus $\Lambda$  is a surjective map.

The part of theorem concerning   $\ker \Lambda \circ \chi $  can be proved in the same manner as  the proof  in the preprint \cite{KW}, see also \cite{HW}.
 \end{proof}
 
 Observe, that $\left({\rm Sym}\left(V_{\mathcal{I,\mathit{\mathbf{d}}}} \right) ^{{\mathfrak{u}_2}}\right)_{\mathit{\mathbf{m}}}=\left(\ker \mathcal{D}_{\mathit{\mathbf{d}}}\right)_{\mathit{\mathbf{m}}},$  where $\mathcal{D}_{\mathit{\mathbf{d}}}$ is a derivation of the polynomial algebra $\mathbb{K}[\mathcal{I}].$
Therefore we  got a handy tool for writing  elements of the kernel of arbitrary linear locally nilpotent derivations.

It is best to look now at some examples.

{\bf Example 3.1.} Let us to consider the symbolic expression  $P:=[\alpha, \beta]^2 [\beta, \gamma] \gamma_0^2.$
Then  ${\rm supp}\, P=\left\{\alpha, \beta, \gamma \right\},$ ${\rm wt}\,P=(2,3,3)=(2)^1 (3)^2,$ ${\rm supp}\, P/\sim=\{\{\alpha\}, \{\beta, \gamma\} \},$ $|{\rm supp}\, P/\sim|=2.$ Put $\mathcal{I}=\{ x, y \}$ and associate $\alpha$  to $x$ and both $\beta, \gamma$ associate to $y.$   Since $\mathit{\mathbf{d}}=(2,3)$ then the map $\Lambda$ acts by
\begin{align*}
&\alpha_0^2 \mapsto x_0,\alpha_0 \alpha_1 \mapsto x_1, \alpha_1^2 \mapsto 2!\,x_2,\\
& \beta_0^3  \mapsto y_0, \beta_0^2 \beta_1  \mapsto y_1,\beta_0 \beta_1^2  \mapsto 2!\,y_2, \beta_1^3  \mapsto 3!\,y_3, \\
& \gamma_0^3  \mapsto y_0, \gamma_0^2 \gamma_1  \mapsto y_1,\gamma_0 \gamma_1^2  \mapsto 2!\,y_2, \gamma_1^3  \mapsto 3!\,y_3.
\end{align*}
We  have 
\begin{gather*}
P:=[\alpha, \beta]^2 [\beta, \gamma] \gamma_0^2={\alpha_{{0}}}^{2}{\beta_{{1}}}^{2}\beta_{{0}}{\gamma_{{0}}}^{2}\gamma
_{{1}}-{\alpha_{{0}}}^{2}{\beta_{{1}}}^{3}{\gamma_{{0}}}^{3}-2\,\alpha_{{0}}\alpha_{{1}}\beta_{{1}}{\beta_{{0}}}^{2}
{
\gamma_{{0}}}^{2}\gamma_{{1}}+\\+2\,\alpha_{{0}}\alpha_{{1}}{\beta_{{1}}}^{2}\beta_{
{0}}{\gamma_{{0}}}^{3}+{\alpha_{{1}}}^{2}{\beta_{{0}}}^{3}
{\gamma_{{0}}}^{2}\gamma_{{1}}-{\alpha_{{1}}}^{2}{\beta_{{0}}}^{2}
\beta_{{1}}{\gamma_{{0}}}^{3}.
\end{gather*}
Thus
\begin{gather*}
\Lambda(P)=2\,x_0 y_2 y_1-6\,x_0 y_3 y_0-2\,x_1 y_1^2
+4\,x_1 y_2 y_0+2\,x_2 y_0 y_1
-2\,x_2 y_1y_0=\\
=2\,x_0 y_2 y_1-6\,x_0 y_3 y_0-2\,x_1 y_1^2
+4\,x_1 y_2 y_0.
\end{gather*}
Consider the polynomial algebra  $\mathbb{K}[X_2,Y_3]:=\mathbb{K}[x_0,x_1,x_2,y_1,y_2,y_3].$ The polynomial  $\Lambda(P)$  belongs to the kernel of the derivation $\mathcal{D}_{(2,3)}$ of the algebra $\mathbb{K}[X_2,Y_3]$ defined by 
\begin{align*}
&\mathcal{D}_{(2,3)}(x_i)=x_{i-1}, i=0,1,2,\mathcal{D}_{(2,3)}(y_i)=y_{i-1}, i=0,1,2,3, \mathcal{D}_{(2,3)}(x_0)=\mathcal{D}_{(2,3)}(y_0)=0.
\end{align*}

{\bf Example 3.2.} Let $P=[\alpha, \beta]^n.$ Then  ${\rm supp}\, P=\left\{\alpha, \beta \right\},$ ${\rm wt}\,P=(n,n)=(n)^2,$ ${\rm supp}\, P/\sim=\{\{\alpha, \beta \}\},$ $|{\rm supp}\, P/\sim|=1.$ Put $\mathcal{I}=\{ x \}$ and associate $\alpha, \beta$  to $x.$ The map $\Lambda$  acts by 
$$
\alpha_0^{n-i}\alpha_1^i \mapsto i!\,x_i, \beta_0^{n-i}\beta_1^i \mapsto i!\,x_i, i=0,1,\ldots,n.
$$
We have 
$$
\Lambda\left([\alpha, \beta]^n\right)=\Lambda\left(\sum_{i=0}^n (-1)^i \alpha_0^{n-i}\alpha_1^i \beta_0^{i}\beta_1^{n-i}\right)=\sum_{i=0}^n (-1)^i i! \,(n-i)!\, x_i x_{n-i}.
$$
Observe that $[\alpha,\beta]^n=(-1)^n [\beta, \alpha]^n$   but, obviously, $\Lambda\left([\alpha, \beta]^n\right)=\Lambda\left([\beta, \alpha]^n\right).$ It follows that $\Lambda\left([\alpha, \beta]^n\right)=0,$  for odd $n.$ Note, that the polynomial belongs to the kernel of the basic Weitzenb\"ock derivation $\mathcal{D}_{n}$ of $\mathbb{K}[V_{x,n}]$ defined by 
$\mathcal{D}_{n}(x_i)=x_{i-1}.$

{\bf Example 3.3.} Consider the polynomial $A=3\,x_{{1}}x_{{2}}x_{{0}}-3\,x_{{3}}{x_{{0}}}^{2}-{x_{{1}}}^{3} \in \ker \mathcal{D}_{3}.$  Find a symbolic expression of the semi-invariant $A.$ To get multilinear polynomial  polarize $A$ two times with respect to letters $y$  and $z:$
\begin{align*}
&P_y(A)= 3\,y_{{0}}x_{{1}}x_{{2}}-6\,y_{{0}}x_{{3}}x_{{0}}+3\,y_{{1}}x_{{0}}x_{
{2}}-3\,y_{{1}}{x_{{1}}}^{2}+3\,y_{{2}}x_{{0}}x_{{1}}-3\,y_{{3}}{x_{{0
}}}^{2},\\
&\mathcal{P}(A)=P_z(P_y(A))=3\,z_{{0}}y_{{1}}x_{{2}}-6\,z_{{0}}y_{{0}}x_{{3}}+3\,z_{{0}}y_{{2}}x_
{{1}}-6\,z_{{0}}y_{{3}}x_{{0}}+3\,z_{{1}}y_{{0}}x_{{2}}-6\,z_{{1}}y_{{
1}}x_{{1}}+\\&+3\,z_{{1}}y_{{2}}x_{{0}}+3\,z_{{2}}y_{{0}}x_{{1}}+3\,z_{{2}
}y_{{1}}x_{{0}}-6\,z_{{3}}y_{{0}}x_{{0}}
.
\end{align*}
The polynomial $\mathcal{P}(A)$ has the multidegree $(1,1,1).$
The map $\widetilde{\Lambda}$ acts  by 
$$
x_i \mapsto 1/i!\,\alpha_0^{3-i} \alpha_1^i, y_i \mapsto 1/i!\,\beta_0^{3-i} \beta_1^i, z_i \mapsto 1/i!\,\gamma_0^{3-i} \gamma_1^i.
$$
We have 
\begin{gather*}
\widetilde{\Lambda}\left(\mathcal{P}(A)\right)=-{\gamma_{{0}}}^{3}{\beta_{{0}}}^{3}{\alpha_{{1}}}^{3}+3/2\,{\gamma_{{0
}}}^{3}{\beta_{{0}}}^{2}\beta_{{1}}\alpha_{{0}}{\alpha_{{1}}}^{2}+3/2
\,{\gamma_{{0}}}^{3}\beta_{{0}}{\beta_{{1}}}^{2}{\alpha_{{0}}}^{2}
\alpha_{{1}}-{\gamma_{{0}}}^{3}{\beta_{{1}}}^{3}{\alpha_{{0}}}^{3}+\\+3/2
\,{\gamma_{{0}}}^{2}\gamma_{{1}}{\beta_{{0}}}^{3}\alpha_{{0}}{\alpha_{
{1}}}^{2}-6\,{\gamma_{{0}}}^{2}\gamma_{{1}}{\beta_{{0}}}^{2}\beta_{{1}
}{\alpha_{{0}}}^{2}\alpha_{{1}}+3/2\,{\gamma_{{0}}}^{2}\gamma_{{1}}
\beta_{{0}}{\beta_{{1}}}^{2}{\alpha_{{0}}}^{3}+3/2\,\gamma_{{0}}{
\gamma_{{1}}}^{2}{\beta_{{0}}}^{3}{\alpha_{{0}}}^{2}\alpha_{{1}}+\\+3/2\,
\gamma_{{0}}{\gamma_{{1}}}^{2}{\beta_{{0}}}^{2}\beta_{{1}}{\alpha_{{0}
}}^{3}-{\gamma_{{1}}}^{3}{\beta_{{0}}}^{3}{\alpha_{{0}}}^{3}
.
\end{gather*}
After simplification we obtain 
\begin{gather*}
2\,\widetilde{\Lambda}\left(\mathcal{P}(A)\right)=3\,\beta_0 {\gamma_{
{0}}}^{2}{[\alpha, \beta]}^{2}[\alpha, \gamma]+3\,{\beta_0}^{2}\gamma_{{0}}[\alpha, \beta]{[\alpha, \gamma]}^{2}-2\,{\beta_0}^{3}{[\alpha, \gamma]}^{3}-2\,{\gamma_{{0}}}^{3}{[\alpha, \beta]}^{3}.
\end{gather*}
Taking into account $\Lambda({{\gamma_{{0}}}^{3}[\alpha, \beta]}^{3})=\Lambda({\beta_0}^{3}{[\alpha, \gamma]}^{3})=0,$ and $$\Lambda(\beta_0 {\gamma_{
{0}}}^{2}{[\alpha, \beta]}^{2}[\alpha, \gamma])=\Lambda({\beta_0}^{2}\gamma_{{0}}[\alpha, \beta]{[\alpha, \gamma]}^{2}),$$ we get that $\widetilde{\Lambda}\left(\mathcal{P}(A)\right)=3\,{\beta_0 {\gamma_{
{0}}}^{2} [\alpha, \beta]}^{2}[\alpha, \gamma].$ Thus ${\beta_0 {\gamma_{
{0}}}^{2} [\alpha, \beta]}^{2}[\alpha, \gamma]$ is a symbolic expression of $A$ and 
$$
A=\Lambda \left({\beta_0 {\gamma_{
{0}}}^{2} [\alpha, \beta]}^{2}[\alpha, \gamma]\right).
$$

{\bf Example 3.4.} Put $P=\prod_{\alpha < \beta}[\alpha,\beta]^2,$  $|\,{\rm supp} P|=2\,k.$ We have, see proof in  \cite{KR}, that 
$$
\Lambda\left(\prod_{\alpha < \beta}[\alpha,\beta]^2 \right)=(k+1)! \begin{vmatrix} x_0 &  x_1 & 2!\, x_2 &\cdots  & k!\, x_k  \\
x_1 & 2 x_2 &3! x_3 &\cdots  & (k+1)! x_{k+1}  \\
\hdotsfor{5}\\
 (k-1)!x_{k-1} &  k!x_{k} & (k+1)!x_{k+1} &\cdots  &(2k-1)!x_{2k-1}\\
 k! x_k &  (k+1)!\,x_{k+1} & (k+2)!\,x_{k+2} &\cdots  & (2k)!\,x_{2k}  
\end{vmatrix}.
$$
This  semi-invariant belongs to $\ker \mathcal{D}_{2\,k},$ has degree $k+1$ and  called the {\it  catalecticant}.

%%%%%%%%%%%%%%%%%%%%%%%%%%%%%%%%%%%%%555
\section{Convolution and semi-tranvectant}
%%%%%%%%%%%%%%%%%%%%%%%%%%%%%%%%%%%%%%%%

There are simple and effective way to find  semi-invariants  of given multidegree ${\bf m.}$ 
The following differential  operator on    ${\rm Sym}\left(V_{\mathcal{J}} \right):$  
 $$
{\rm Conv}_{\alpha, \beta}:= [\alpha,\beta] \frac{\partial^2}{\partial \alpha_0\partial \beta_0}, \alpha, \beta \in \mathcal{J}
$$
is called the  { \textit{ convolution  }} with respect to the symbol letters $\alpha$ and  $\beta.$
Obviously, the convolution operator  does not change the weight of a symbolic  expression, so
  ${\rm Conv}_{\alpha, \beta}$ is an endomorphism of the vector space   ${\rm Sym}\left(V_{\mathcal{J}}^{{\mathfrak{u}_2}} \right)_{{\bf d}^{{\bf m}}}.$ 

{\bf Example 4.1} Let us  find an elements of kernel of the derivation  $\mathcal{D}_{1,2,3}$  of the multidegree  $(1,2,1)$ and the weight $(1,2,2,3)=(1)^1 (2)^2 (3)^1.$ Put $\mathcal{J}=\{\alpha, \beta, \gamma,\delta \},$ $\mathcal{I}=\{x,y, z\}.$ To the symbol letter  $\alpha$ we  associate  the roman letter $x,$ to both letters  $\beta, \gamma$ we associate the letter  $y$ and  to symbol letter  $\delta$ we  associate  the letter $z.$
 The  map $\Lambda$ acts  by  $\alpha_i \mapsto x_i,$ $i=0,1,$ $\beta_0^{2-i} \beta_1^i \mapsto i!\, y_i,$
$\gamma_0^{2-i} \gamma_1^i \mapsto i!\, y_i,$ $i=0,1,2$  and $\delta_0^{3-i} \delta_1^i \mapsto i!\, z_i,$ $z=0,1,2,3.$
It is clear that  the following  symbolic expression  $\Phi=\alpha_0 \beta_0^2 \gamma_0^2 \delta_0^3$  has the weight $(1)^1 (2)^2 (3)^1.$ The polynomial $\Lambda(\Phi)$ is equal to  $x_0 y_0^2 z_0$ and has the multidegree $(1,2,1).$  By direct calculations we get 
\begin{align*}
&{\rm Conv}_{\alpha, \beta}(\Phi)=2\,[\alpha,\beta] \beta_0 \gamma_0^2 \delta_0^3 \longmapsto 2\,y_0 z_0 \begin{vmatrix} x_0 & x_1 \\ y_0 & y_1 \end{vmatrix},\\
&{\rm Conv}_{\gamma, \delta}\left( {\rm Conv}_{\alpha, \beta}(\Phi)\right)=12\,[\alpha,\beta] [\gamma,\delta] \beta_0 \gamma_0 \delta_0^2 \mapsto 12 \begin{vmatrix} x_0 & x_1 \\ y_0 & y_1 \end{vmatrix} \begin{vmatrix} y_0 & y_1 \\ z_0 & z_1 \end{vmatrix},\\
&{\rm Conv}_{\gamma, \delta}^2 \left( {\rm Conv}_{\alpha, \beta}(\Phi)\right)=24\,[\alpha,\beta] [\gamma,\delta]^2 \beta_0  \delta_0 \mapsto 48\,(z_{{0}}y_{{2}}-y_{{1}}z_{{1}}+y_{{0}}z_{{2}})\begin{vmatrix} x_0 & x_1 \\ y_0 & y_1 \end{vmatrix},\\
&{\rm Conv}_{\beta, \delta}\left({\rm Conv}_{\gamma, \delta}^2 \left( {\rm Conv}_{\alpha, \beta}(\Phi)\right)\right)=24\,[\alpha,\beta] [\gamma,\delta]^2 [\beta  \delta] \mapsto 48 \,(3\,x_{{0}}y_{{0}}z_{{3}}y_{{1}}-2\,x_{{0}}y_{{2}}y_{{0}}z_{{2}}-\\
&-2\,x_{
{0}}{y_{{1}}}^{2}z_{{2}}+3\,x_{{0}}y_{{2}}y_{{1}}z_{{1}}-2\,x_{{0}}{y_
{{2}}}^{2}z_{{0}}-3\,{y_{{0}}}^{2}x_{{1}}z_{{3}}+3\,y_{{1}}x_{{1}}y_{{0
}}z_{{2}}-{y_{{1}}}^{2}x_{{1}}z_{{1}}-y_{{0}}x_{{1}}z_{{1}}y_{{2}}+y_{
{1}}x_{{1}}z_{{0}}y_{{2}}).
\end{align*}

For a symbolic expression $P$  denote by  ${\rm Conv}(P)$  the set of its all possible convolutions. For a subalgebra  $\Delta \in {\rm Sym}\left(V_{\mathcal{J}} \right)^{{\mathfrak{u}_2}}$  denote by ${\rm Conv}(\Delta )$ the subalgebra generated by all possible convolutions ${\rm Conv}(P),$ $P \in \Delta.$
The following statement  holds:
\begin{lm}\label{conv}
$$
{\rm Sym}\left(V_{\mathcal{J}} \right)^{{\mathfrak{u}_2}}={\rm Conv}({\rm Sym}\left(\oplus_{\alpha \in \mathcal{J}}\mathbb{K} \alpha_0 \right)).
$$
\end{lm}
\begin{proof}
Let $P \in ({\rm Sym}\left(V_{\mathcal{J}} \right)^{{\mathfrak{u}_2}})_{\bf w}.$  Then $P$ is obtained by convolutions of the semi-invariants $\prod_{\alpha \in \mathcal{J}} \alpha_0^{{\rm wt}_{\alpha }}  \in {\rm Sym}\left(\oplus_{\alpha \in \mathcal{J}}\mathbb{K} \alpha_0 \right).$
\end{proof}

{\bf Example 4.1} Let $\mathcal{J}=\{\alpha, \beta, \gamma \}.$ The component $\left({\rm Sym}\left(V_{\mathcal{J}} \right)^{{\mathfrak{u}_2}}\right)_{(2,1,1)}$  is generated by the following $5$  semi-invariants 
\begin{align*}
 \alpha_0^2 \beta_0 \gamma_0, \alpha_0  \gamma_0 [\alpha, \beta],  \alpha_0  \beta_0 [\alpha, \gamma],  \alpha_0^2 [\beta, \gamma],[\alpha,\gamma] [\alpha, \beta].
\end{align*}
All of them are the convolutions of the semi-invariant  $\alpha_0^2 \beta_0 \gamma_0.$

 Let $F \in {\rm Sym}\left(  V_{\mathcal{I,\mathit{\mathbf{d}}}} \right) ^{{\mathfrak{u}_2}}$  and $\Phi \in {\rm Sym}\left(V_{\mathcal{J}} \right)^{{\mathfrak{u}_2}}$ is its symbolic representation. Denote by ${\rm Conv}_{\Lambda}(F)$  the set of elements $\Lambda \left({\rm Conv}(\Phi) \right).$ The elements of ${\rm Conv}_{\Lambda}(F)$  are called the $\Lambda$-{\it convolutions.} For a subalgebra  $\mathcal{T} \subset {\rm Sym}\left(  V_{\mathcal{I,\mathit{\mathbf{d}}}} \right) ^{{\mathfrak{u}_2}} $  denote by  ${\rm Conv}_{\Lambda}(\mathcal{T})$ the algebra  generated by all $\Lambda$-convolutions of all elements of the algebra $\mathcal{T}.$

\begin{te}\label{sym_kryt}   Let  $\mathcal{N}_{\mathcal{I}}:=\mathbb{K}[x_0 \mid x \in \mathcal{I}].$ Then   ${\rm Sym}\left(  V_{\mathcal{I,\mathit{\mathbf{d}}}} \right) ^{{\mathfrak{u}_2}}={\rm Conv}_{\Lambda}(\mathcal{N}_{\mathcal{I}}).$
\end{te}
\begin{proof}
  Consider an arbitrary homogeneous semi-invariant  $F \in {\rm Sym}\left(  V_{\mathcal{I,\mathit{\mathbf{d}}}} \right) ^{{\mathfrak{u}_2}}.$ Let  $\Phi$  be its symbolic expression.   By Lemma  \ref{conv},  $\Phi$  belongs to   ${\rm Conv}({\rm Sym}\left(\oplus_{\alpha \in \mathcal{J}}\mathbb{K} \alpha_0 \right))$ and $|\mathcal{J}|=\deg F.$ It is easy to see, that   $$\Lambda({\rm Sym}\left(\oplus_{\alpha \in \mathcal{J}}\mathbb{K} \alpha_0 \right))=\mathcal{N}_{\mathcal{I}}   \text{ and } \Lambda \left(  {\rm Conv}({\rm Sym}\left(\oplus_{\alpha \in \mathcal{J}}\mathbb{K} \alpha_0 \right)) \right)={\rm Conv}_{\Lambda}(\mathcal{N}_{\mathcal{I}}),$$  thus $\Lambda(\Phi)=F \in {\rm Conv}(\mathcal{N}_{\mathcal{I}}).$  We  get  ${\rm Sym}\left(  V_{\mathcal{I,\mathit{\mathbf{d}}}} \right) ^{{\mathfrak{u}_2}} \subseteq {\rm Conv}_{\Lambda}\left(\mathcal{N}_{\mathcal{I}}\right).$  The inclusion ${\rm Conv}_{\Lambda}\left(\mathcal{N}_{\mathcal{I}}\right) \subseteq {\rm Sym}\left(  V_{\mathcal{I,\mathit{\mathbf{d}}}} \right) ^{{\mathfrak{u}_2}} $   is obvious.
\end{proof}
We  have got a tool  to find a generating system of kernel for  Weitzenb\"ock derivations.

{\bf Example 4.2}  Let $|\mathcal{I}|=n$ and ${\mathbf d} = (1,1,\ldots, 1).$  Prove that   $$\ker \mathcal{D}_{\mathit{\mathbf{d}}}=\mathbb{K}[x_0, x_0 y_1-x_1 y_0 \mid x \neq y, x,y \in \mathcal{I}].$$
In fact, let $F$ be a homogeneous polynomial of  $\mathcal{N}_{\mathcal{I}}$ of degree $m.$ Then its symbolic representation  $\Phi$ has the form $\Phi=\prod_{\alpha \in \mathcal{J}} \alpha_0,$ $\mathcal{J} \subset \mathfrak{G},$ $ |\mathcal{J}|=m.$
Since all factors of $\Phi$ have the degrees 1,  we obtain that all possible convolutions of the polynomial  $\Phi$  have form $\prod_{\alpha \in \mathcal{J}} \alpha_0 \prod_{\beta, \gamma \in \mathcal{J}}[\beta, \gamma],$ for $\alpha \neq \beta , \gamma.$ It follows
$$
\Lambda \left(\prod_{\alpha \in \mathcal{J}} \alpha_0 \prod_{\beta, \gamma \in \mathcal{J}}[\beta, \gamma] \right)=\prod_{x \in \mathcal{I}} x_0 \prod_{y,z \in \mathcal{I}}(y_0 z_1 -y_1 z_0), x \neq y,z.
$$
Thus $ {\rm Conv}_{\Lambda}(F) \in \mathbb{K}[x_0, x_0 y_1-x_1 y_0 \mid x \neq y, x,y \in \mathcal{I}]$ and $${\rm Conv}_{\Lambda}(\mathcal{N}_{\mathcal{I}})=\mathbb{K}[x_0, x_0 y_1-x_1 y_0 \mid x \neq y, x,y \in \mathcal{I}].$$

{\bf Example 4.3}  Let  $|\mathcal{I}|=n \geq 3,$ ${\mathbf d} = (2,2,\ldots, 2)$ and let  $F$ be an homogeneous  polynomial of  $\mathcal{N}_{\mathcal{I}}$ of degree $m.$ Then its symbolic representation  $\Phi$ has the form $\Phi=\prod_{\alpha \in \mathcal{J}} \alpha_0^2,$ $\mathcal{J} \subset \mathfrak{G},$ $ |\mathcal{J}|=m \cdot n.$ 
For  $m=1$ we have  $\Phi=\alpha_0^2, \alpha \in \mathcal{J}.$  It is obvious that ${\rm Conv}(\Phi)=\{ \alpha_0^2, \alpha \in \mathcal{J}\}.$ 

For  $m=2$ we  have  $\Phi=\alpha_0^2 \beta_0^2, \alpha \neq \beta, \alpha, \beta \in \mathcal{J}.$  There are only two convolutions of $\Phi:$ 
$\alpha_0 \beta_0 [\alpha, \beta]$ and $[\alpha, \beta]^2.$

For  $m=3$ we have $\Phi=\alpha_0^2 \beta_0^2 \gamma_0^2, \alpha, \beta, \gamma \in \mathcal{J}.$  There exists a  unique  not decomposable convolution: ${[\alpha, \beta] [\alpha, \gamma] [\beta, \gamma].}$

 If  $m>3$ then any symbol expression is either decomposable or belongs to $\ker \Lambda,$ see \cite{GrY}, page.162.
We  have 
\begin{gather*}
\Lambda(\alpha_0^2)=x_0, \Lambda \left(\alpha_0 \beta_0 [\alpha, \beta] \right)=\begin{vmatrix} x_0 & x_1 \\ y_0 & y_1 \end{vmatrix},\Lambda \left([\alpha, \beta]^2 \right)=2\,(x_0 y_2-x_1 y_1+ x_2 y_2),
\end{gather*}
\begin{gather*}
\Lambda\left([\alpha, \beta] [\alpha, \gamma] [\beta, \gamma]\right)=\Lambda({\alpha_{{0}}}^{2}\beta_{{1}}{\gamma_{{1}}}^{2}\beta_{{0}}-{\alpha_{{0
}}}^{2}{\beta_{{1}}}^{2}\gamma_{{1}}\gamma_{{0}}+\\+\alpha_{{0}}{\beta_{{
1}}}^{2}{\gamma_{{0}}}^{2}\alpha_{{1}}-{\beta_{{0}}}^{2}\alpha_{{1}}
\alpha_{{0}}{\gamma_{{1}}}^{2}+{\beta_{{0}}}^{2}{\alpha_{{1}}}^{2}
\gamma_{{0}}\gamma_{{1}}-\beta_{{0}}{\alpha_{{1}}}^{2}{\gamma_{{0}}}^{
2}\beta_{{1}})=\\=2\,x_{{0}}z_{{2}}y_{{1}}-2\,x_{{0}}y_{{2}}z_{{1}}+2\,x_{{1}}y_{{2}}z_{
{0}}-2\,y_{{0}}x_{{1}}z_{{2}}+2\,y_{{0}}x_{{2}}z_{{1}}-2\,y_{{1}}x_{{2
}}z_{{0}}
=\\
=2 \,\begin{vmatrix} x_0 & y_0 & z_0 \\ x_1 & y_1 &z_1 \\ x_2 & y_2 & z_2 \end{vmatrix}.
\end{gather*}
Thus, the kernel 
$
\ker \mathcal{D}_{(2,2,\ldots, 2)}
$
is generated by the semi-invariants of  these four types.

With increasing of $d_i$ it becomes difficult to apply the Theorem \ref{sym_kryt}. We will offer other similar but  more effective approach to find the kernel of a Weitzenb\"ock derivations.
In the paper \cite{Bed_In7} we introduced the conception of semi-transvectant, an analogue of the classical transvectant. 
Recal that the algebra  ${\mathfrak{u}_2}$ acts on    $ {\rm Sym}\left(  V_{\mathcal{I,\mathit{\mathbf{d}}}} \right)$ by the locally nilpotent derivation ${\bf D}=\mathcal{D}_{\mathit{\mathbf{d}}}.$
Let  $F,G \in {\rm Sym}\left(  V_{\mathcal{I,\mathit{\mathbf{d}}}} \right) ^{{\mathfrak{u}_2}}$  be two semi-invariants of degrees   $p$ and $q,$ respectively. The semi-invariant of the form   
\begin{equation}\label{trans}
[F,G]^r:= \sum_{i=0}^r (-1)^i { r \choose i } \frac{{\bf D}_*^i(F)}{[p]_i}  \frac{{\bf D}_*^{r-i}(G)}{[q]_{r-i}}, 
\end{equation}
$\mbox{ } 0 \leq r \leq \min(p,q),$
$[m]_i:=m (m-1) \ldots (m-(i-1)), m \in \mathbb{Z},$
 is called the  {\it $r$-th semi-transvectant} of the semi-invariants   $F$ and $G.$

{\bf Example   4.4. } The semi-transvectant  $[F,G]^1:=[F,G]$  is called {\it the semi-Jacobian}. If $F, G, H$ are three semi-invariants of orders  greater than unity, then the iterated semi-Jacobian $[[F, G], H]$ is reducible \cite{Gle}  and 
$$
[[F, G], H]=\frac{{\rm ord}(F)-{\rm ord}(G)}{2\,({\rm ord}(F)+{\rm ord}(G)-2)}+\frac{1}{2} [F,G]^2 H+\frac{1}{2} [F,H]^2 G-\frac{1}{2} [G,H]^2 F.
$$

 {\bf Example   4.5. } The semi-invariant $[F,F]^2:={\rm Hes}(F)$  is called {\it the semi-Hessian}. The square of a semi-Jacobian $[F,G]$ is given by the formula
 $$
 [F,G][F,G]=[F,G]^2 F\,G-\frac{1}{2} {\rm Hes}(F) G^2 -\frac{1}{2} {\rm Hes}(G) F^2.
 $$
 
 {\bf Example   4.6. } We  have
 \begin{gather*}
 \left( \alpha_0^2,\beta_0 \gamma_0^2  \right)^2=\frac{1}{3}\,{\alpha_{{0}}}^{2}\beta_{{1}}\gamma_{{1}}\gamma_{{0}}+\frac{1}{6}\,{
\alpha_{{0}}}^{2}\beta_{{0}}{\gamma_{{1}}}^{2}-\frac{1}{3}\,\alpha_{{1}}\alpha
_{{0}}\beta_{{1}}{\gamma_{{0}}}^{2}-\frac{2}{3}\,\alpha_{{1}}\alpha_{{0}}
\gamma_{{1}}\beta_{{0}}\gamma_{{0}}+\frac{1}{2}\,{\alpha_{{1}}}^{2}\beta_{{0}}
{\gamma_{{0}}}^{2}
=\\
=\frac{1}{2}\,{[\alpha,\gamma]}^{2}\beta_{{0}}-\frac{1}{3}\,\alpha_{{0}}[\alpha,\gamma][\beta,\gamma].
 \end{gather*}

The following statement holds:
\begin{lm}[\cite{Bed_In7}]\label{prop-t} 
\begin{align*}
& \text{({\it i}) for   $ 0 \leq r \leq \min({\rm ord}(x_0), \max( {\rm ord}(F),{\rm ord}(G))$ the semi-transvectant  $[x_0, F\,G]^r$ is reducible  ;}\\
&\text{ ({\it ii}) if  $ {\rm ord}(F)=0,$ \text{  then, }  $[x_0, F\,G]^r=F [x_0,G]^r;$}\\
&\text{({\it iii})}\, \, {\rm ord}([F,G]^r)={\rm ord}(F)+{\rm ord}(G)-2\, r. 
\end{align*}
\end{lm}

There is a close relationship between the convolutions and the semi-transvectants. A $k$-{\it fold contraction}  of two disjoint symbolic
expressions $\Phi$ and $\Psi$ is a symbolic expression $$\left(\prod_{\alpha, \beta} {\rm Conv }_{\alpha, \beta}\right) \left( \Phi\cdot \Psi \right)$$ where the product runs over
$k$ pairs $(\alpha,\beta) \in {\rm Supp}\, \Phi \times {\rm Supp}\, \Psi.$ 

\begin{lm}[\cite{KW}]
 Let $\Phi$ and $\Psi$ be two disjoint symbolic expressions. The semi-transvectant
$[\Phi, \Psi]^k$ is a linear combination of  semi-invariants $T$ where $T$ runs through the k-fold
contractions of $\Phi$ and $\Psi,$  and each such $T$ occurs with a positive rational coefficient
$q_T$ where $\sum_{T}q_T = 1.$
\end{lm}
Using the semi-transvectants we offer an algorithm  for computation of the kernel of Weitzenb\"ock  derivations.
For  $x \in \mathcal{I}$ put   $\tau_{x,i}(F):=(x_0, F)^i,$ ${ i \leq \min \left({\rm ord}(x_0), {\rm ord}(F)\right).}$ 
For a subalgebra   $ T \subseteq {\rm Sym}\left(  V_{\mathcal{I,\mathit{\mathbf{d}}}} \right) ^{{\mathfrak{u}_2}}$  denote by  $\tau(T)$ the algebra generated by the elements   $\tau_{x,i}(F), F\in T,$   $ i \leqslant \min({\rm ord}(x_0),\mbox{ord}(F)).$
The following theorem is a weak form of Gordan's theorem: 
\begin{te}[\cite{B_UMZH}]
Let  $\mathcal{T}$ be a subalgebra of   $\ker \mathcal{D}_{\mathit{\mathbf{d}}}$   containing   $\mathcal{N}_{I}$  and $\tau (\mathcal{T}) \subseteq \mathcal{T}$. Then  $\ker \mathcal{D}_{\mathit{\mathbf{d}}}=\mathcal{T}.$
\end{te}
The theorem implies the following algorithm for  $\ker \mathcal{D}_{\mathit{\mathbf{d}}}.$  Define  the series of subalgebras 
$$
\mathcal{T}_1 \subseteq \mathcal{T}_2 \subseteq \mathcal{T}_3 \subseteq \cdots \subset \mathcal{T}_k \subseteq \cdots , 
$$
where $\mathcal{T}_1=\mathcal{N}_{\mathcal{I}}$ i $\mathcal{T}_i:=\tau(\mathcal{T}_{i-1}).$
If for  some  $i$ we have that  $\mathcal{T}_i=\mathcal{T}_{i+1},$ then $\mathcal{T}_i=\ker \mathcal{D}_{\mathit{\mathbf{d}}}.$

{\bf Example 4.7} Consider the basic Weitzenb\"ock derivation   $ \mathcal{D}_3,$ $\mathcal{I}=\{ x \}.$ We have  $\mathcal{T}_1=\mathbb{K}[x_0].$ The subalgebra   $\mathcal{T}_2$  generated by the elements $\tau_i(x_0^j).$ By Lemma  \ref{prop-t}, $(i)$ for  $j>1$ all of them are reducible one. The only irreducible semi-invariant is 
 ${dv=\tau_2(x_0)}$.  Therefore, $\mathcal{T}_2=\mathbb{K}[x_0, dv].$  The algebra  $\mathcal{T}_3$ generated by  $\tau_i( x_0^k dv^l),$ $i,k,l \in \mathbb{N}.$ Since ${{\rm ord}({dv})=2}$, then Lemma \ref{prop-t},(i) implies that   the algebra $\mathcal{T}_3$  consists only the following elements $\tau_1(dv),\tau_2(dv)$ i $\tau_3(dv^2)$.  By (\ref{trans})  we obtain that   $\tau_3(dv^2)=0,$  $\tau_2(dv)=0 $   and   ${tr:= \tau_1(dv)}\neq 0.$ The direct calculation shows that $tr$  does  not belong to   $T_2$, thus  
$T_3=\mathbb{K}[t,dv,tr]$. The algebra $\mathcal{T}_4$ generated by  $\tau_i( x_0^k dv^l tr^m),$ $i,k,l,m \in \mathbb{N}.$ As above we find the only new element ${ch=\tau_3(tr)}$  for $\mathcal{T}_4.$  We  have   ${{\rm ord}(ch)=0}$  but the algebra  $\mathcal{T}_3$ does not consist any  invariants. It implies that   $\mathcal{T}_4=\mathbb{K}[t,dv,tr,ch]$. 
By lemma  \ref{prop-t}, $(ii)$  we have that the algebra  $\mathcal{T}_5$  does not  consist of any  new semi-invariants. Thus  $\mathcal{T}_5=T_4$   and  $\ker {\bf D}=\mathbb{K}[t,dv,tr,ch]$, where 
$$
\begin{array}{l}
dv=  {x_{1}}^{2}- 2\,x_0\,{x_{2}}, \\
tr= 3\,x_{{3}}{x_{{0}}}^{2}+{x_{{1}}}^{3}-3\,x_{{0}}x_{{1}}x_{{2}},\\
ch=8\,x_{{0}}{x_{{2}}}^{3}+9\,{x_{{3}}}
^{2}{x_{{0}}}^{2}+6\,{x_{{1}}}^{3}x_{{3}}-3\,{x_{{1}}}^{2}{x_{{2}}}^{2
}-18\,x_{{0}}x_{{1}}x_{{2}}x_{{3}}.
\end{array}
$$
Up  to a  constant factor the symbolic representation of  these semi-invariants have the form  $\alpha_0 \beta_0 [\alpha, \beta]^2,$ ${\beta_0 {\gamma_{
{0}}}^{2} [\alpha, \beta]}^{2}[\alpha, \gamma]$ and  $[\alpha,\beta]^2 [\alpha, \gamma] [\beta, \delta][\gamma, \delta]^2,$  respectively.

For    $d=4,5,6$ the kernel of the basic Weitzenb\"ock derivation was calculated in \cite{B_UMZH}. The cases   $d=7,8$ considered in   \cite{Bed_Cov7}, \cite{Bed_Cov8}. For $d>8$ the problem is still  open however  the corresponding algebras of invariants was calculated for  $d=9,10$ in  \cite{Br9}, \cite{Br10}.

{\bf Example 4.8} Consider the derivation  $ \mathcal{D}_{(2,3)},$ $\mathcal{I}=\{x,y\},$ ${\rm ord}(x_0)=2,$ ${\rm ord}(y_0)=3.$ To the set  $\mathcal{J}_1:=\{\alpha, \beta, \gamma, \delta  \}$ we associate the letter $x$ and to set í³ $\mathcal{J}_2:=\{ \varepsilon,\varkappa,\eta,\mu \}$ we associate the letter $y.$  The maps $\Lambda$  acts  by 
$\nu_0^{2-i} \nu_1^i \mapsto i!\, x_i,$ $i=0,1,2$  for $\nu \in \mathcal{J}_1$ and by $\nu_0^{3-i} \nu_1^i \mapsto i!\, y_i,$ $i=0,1,2,3$  for $\nu \in \mathcal{J}_2.$  Put $D:=\Lambda \left([\alpha,\beta]^2 \right),$ $\Delta:=\Lambda \left([\epsilon,\varkappa]^2\epsilon_0 \varkappa_0 \right),$ ${Q:=\Lambda \left([\varepsilon,\varkappa]^2 [\eta, \varepsilon]\varkappa_0 \eta_0^3 \right),}$ $R:=\Lambda \left([\varepsilon,\varkappa]^2 [\varepsilon, \eta] [\varkappa, \mu][\eta, \mu]^2 \right).$

The minimal generating set for   $\ker  \mathcal{D}_{(2,3)}$ consists of the following  $15$ elements:
\begin{gather*}
x_0,y_0,\\
D, [x_0,y_0]^2, \Delta, [x_0,y_0],\\
[x_0, \Delta]^2,[x_0^2,y_0]^3, [x_0,\Delta], Q,\\
R, [x_0,Q]^2,\\
[x_0^3,y_0^2]^6,[x_0^2,Q]^3,\\
[x_0^3,y_0 Q]^6.
\end{gather*}
The proof for the corresponding algebras of covariants see  in \cite{Gor} or in \cite{GrY},\cite{Gle}.

For instance, let us calculate the explicit form of the semi-invariant  $(x_0,Q)^2.$ We  have
\begin{gather*}
\left( \alpha_0^2 , [\varepsilon,\varkappa]^2 [\eta, \varepsilon]\varkappa_0 \eta_0^3  \right)^2 = [\varepsilon,\varkappa]^2 [\eta, \varepsilon] \left( \alpha_0^2,\varkappa_0 \eta_0^3  \right)^2 =
\end{gather*}
\begin{multline*}
=[\varepsilon,\varkappa]^2 [\eta, \varepsilon] (\frac{1}{6} \,{\alpha_{{0}}}^{2}\varkappa_{{0}}{\eta_{{1}}}^{2}+\frac{1}{3}\,{\alpha_{{0}}
}^{2}\varkappa_{{1}}\eta_{{1}}\eta_{{0}}-\frac{2}{3}\,\alpha_{{1}}\alpha_{{0}}
\eta_{{1}}\varkappa_{{0}}\eta_{{0}}-\frac{1}{3}\,\alpha_{{1}}\alpha_{{0}}\varkappa_{
{1}}{\eta_{{0}}}^{2}+\frac{1}{2}\,\varkappa_{{0}}{\eta_{{0}}}^{2}{\alpha_{{1}}}^{
2}
)= 
\end{multline*}
\begin{align*}
&=\frac{5}{6}\,{\alpha_{{0}}}^{2}\varkappa_{{0}}{\eta_{{1}}}^{2}{\varepsilon_{{0}}}^{2
}{\varkappa_{{1}}}^{2}\eta_{{0}}\varepsilon_{{1}}-\frac{2}{3}\,{\alpha_{{0}}}^{2}{
\varkappa_{{0}}}^{2}{\eta_{{1}}}^{2}\varepsilon_{{0}}\varkappa_{{1}}{\varepsilon_{
{1}}}^{2}\eta_{{0}}-\frac{2}{3}\,{\alpha_{{0}}}^{2}{\varkappa_{{1}}}^{2}\eta_{{1}
}{\eta_{{0}}}^{2}\varepsilon_{{0}}\varkappa_{{0}}{\varepsilon_{{1}}}^{2}+\\
&+\frac{2}{3}\,
\alpha_{{1}}\alpha_{{0}}{\eta_{{1}}}^{2}\varkappa_{{0}}\eta_{{0}}{
\varepsilon_{{0}}}^{3}{\varkappa_{{1}}}^{2}+\frac{2}{3}\,\alpha_{{1}}\alpha_{{0}}{
\eta_{{1}}}^{2}{\varkappa_{{0}}}^{3}\eta_{{0}}{\varepsilon_{{1}}}^{2}
\varepsilon_{{0}}+\frac{2}{3}\,\alpha_{{1}}\alpha_{{0}}{\varkappa_{{1}}}^{2}{\eta_{{0
}}}^{3}\varepsilon_{{0}}\varkappa_{{0}}{\varepsilon_{{1}}}^{2}+\\&+
{\varkappa_{{0}}}^{
2}{\eta_{{0}}}^{2}{\alpha_{{1}}}^{2}{\varepsilon_{{0}}}^{2}\varkappa_{{1}}
\varepsilon_{{1}}\eta_{{1}}+\frac{1}{2}\,{\varkappa_{{0}}}^{3}{\eta_{{0}}}^{3}{
\alpha_{{1}}}^{2}{\varepsilon_{{1}}}^{3}-\frac{1}{6}\,{\alpha_{{0}}}^{2}\varkappa_{{0
}}{\eta_{{1}}}^{3}{\varepsilon_{{0}}}^{3}{\varkappa_{{1}}}^{2}+\\&+
\frac{1}{6}\,{\alpha_
{{0}}}^{2}{\varkappa_{{0}}}^{3}{\eta_{{1}}}^{2}{\varepsilon_{{1}}}^{3}\eta_{
{0}}-\frac{1}{6}\,{\alpha_{{0}}}^{2}{\varkappa_{{0}}}^{3}{\eta_{{1}}}^{3}{
\varepsilon_{{1}}}^{2}\varepsilon_{{0}}-\frac{1}{3}\,{\alpha_{{0}}}^{2}{\varkappa_{{1}}
}^{3}{\eta_{{1}}}^{2}\eta_{{0}}{\varepsilon_{{0}}}^{3}+\\&+
\frac{1}{3}\,{\alpha_{{0}}}^{2}{\varkappa_{{1}}}^{3}\eta_{{1}}{
\eta_{{0}}}^{2}{\varepsilon_{{0}}}^{2}\varepsilon_{{1}}+\frac{1}{3}\,{\alpha_{{0}}}^
{2}\varkappa_{{1}}\eta_{{1}}{\eta_{{0}}}^{2}{\varkappa_{{0}}}^{2}{\varepsilon_{
{1}}}^{3}-\frac{2}{3}\,\alpha_{{1}}\alpha_{{0}}\eta_{{1}}{\varkappa_{{0}}}^{3}{
\eta_{{0}}}^{2}{\varepsilon_{{1}}}^{3}+\\&+
\frac{1}{3}\,\alpha_{{1}}\alpha_{{0}}{\varkappa_{{1}}}^{3}{\eta_{{0}}}^{2}{\varepsilon
_{{0}}}^{3}\eta_{{1}}-\frac{1}{3}\,\alpha_{{1}}\alpha_{{0}}\varkappa_{{1}}{\eta_{
{0}}}^{3}{\varkappa_{{0}}}^{2}{\varepsilon_{{1}}}^{3}+\frac{1}{2}\,\varkappa_{{0}}{\eta
_{{0}}}^{3}{\alpha_{{1}}}^{2}{\varepsilon_{{0}}}^{2}{\varkappa_{{1}}}^{2}
\varepsilon_{{1}}-\\&-
\frac{1}{2}\,\varkappa_{{0}}{\eta_{{0}}}^{2}{\alpha_{{1}}}^{2}{
\varepsilon_{{0}}}^{3}{\varkappa_{{1}}}^{2}\eta_{{1}}-{\varkappa_{{0}}}^{2}{
\eta_{{0}}}^{3}{\alpha_{{1}}}^{2}\varepsilon_{{0}}\varkappa_{{1}}{\varepsilon_{
{1}}}^{2}-\frac{1}{2}\,{\varkappa_{{0}}}^{3}{\eta_{{0}}}^{2}{\alpha_{{1}}}^{2}{
\varepsilon_{{1}}}^{2}\varepsilon_{{0}}\eta_{{1}}-\\&-
\frac{4}{3}\,\alpha_{{1}}\alpha_{{0
}}\eta_{{1}}\varkappa_{{0}}{\eta_{{0}}}^{2}{\varepsilon_{{0}}}^{2}{\varkappa_{{
1}}}^{2}\varepsilon_{{1}}-\frac{4}{3}\,\alpha_{{1}}\alpha_{{0}}{\eta_{{1}}}^{2}{
\varkappa_{{0}}}^{2}\eta_{{0}}{\varepsilon_{{0}}}^{2}\varkappa_{{1}}\varepsilon_{{
1}}+\frac{5}{3}\,\alpha_{{1}}\alpha_{{0}}\eta_{{1}}{\varkappa_{{0}}}^{2}{\eta_{{0
}}}^{2}\varepsilon_{{0}}\varkappa_{{1}}{\varepsilon_{{1}}}^{2}+\\
&+\frac{1}{3}\,{\alpha_{{0}}
}^{2}{\varkappa_{{0}}}^{2}{\eta_{{1}}}^{3}{\varepsilon_{{0}}}^{2}\varkappa_{{1}
}\varepsilon_{{1}}-\frac{1}{3}\,\alpha_{{1}}\alpha_{{0}}{
\varkappa_{{1}}}^{3}{\eta_{{0}}}^{3}{\varepsilon_{{0}}}^{2}\varepsilon_{{1}}.
\end{align*}
Then 
\begin{gather*}
(x_0,Q)^2 =\Lambda \left(\left( \alpha_0^2 , [\varepsilon,\varkappa]^2 [\eta, \varepsilon]\varkappa_0 \eta_0^3  \right)^2\right)=\\
=6\,x_{{0}}{y_{{1}}}^{2}y_{{3}}-6\,y_{{2}}y_{{0}}x_{{2}}y_{{1}}-6\,x_{{
1}}y_{{0}}y_{{1}}y_{{3}}-2\,x_{{0}}y_{{1}}{y_{{2}}}^{2}+6\,{y_{{0}}}^{
2}x_{{2}}y_{{3}}-6\,x_{{0}}y_{{2}}y_{{3}}y_{{0}}+8\,x_{{1}}y_{{0}}{y_{
{2}}}^{2}-\\-2\,x_{{1}}y_{{2}}{y_{{1}}}^{2}+2\,{y_{{1}}}^{3}x_{{2}}.
\end{gather*}

 In  \cite{Gor}, \cite{GrY}, \cite{Gall_33}  were calculated the algebras joint covariants which are isomorphic to the kernels of the following derivations  $\mathcal{D}_{(2,4)},$ $\mathcal{D}_{(3,3)}$ and $\mathcal{D}_{(3,3,3)}.$

\end{document}